\documentclass[conference]{IEEEtran}
\usepackage{cite}
\ifCLASSINFOpdf
\usepackage[pdftex]{graphicx}
\graphicspath{{./pdf/}{./jpeg/}}
\DeclareGraphicsExtensions{.pdf,.jpeg,.png}
\else
\usepackage[dvips]{graphicx}
\graphicspath{{./eps/}}
\DeclareGraphicsExtensions{.eps}
\fi
\usepackage{amsmath}
\usepackage{amsfonts}
\usepackage{multirow}
\usepackage{hhline}
\usepackage{xcolor,colortbl}
\definecolor{Gray}{gray}{0.85}
\usepackage{array}

\usepackage{stfloats}
\usepackage{url}
\begin{document}
\title{Maximizing Investment Value of Small-Scale PV in a Smart Grid Environment}

\author{\IEEEauthorblockN{Jeremy Every, Li Li, Youguang G. Guo}
\IEEEauthorblockA{Faculty of Engineering and Information Technology\\
University of Technology Sydney\\
Ultimo 2007, Australia \\
Email: jeremy.every@student.uts.edu.au}
\and
\IEEEauthorblockN{David G. Dorrell}
\IEEEauthorblockA{College of Agriculture, Engineering and Science\\
University of KwaZulu-Natal\\Durban 4041, South Africa\\
Email: dorrelld@ukzn.ac.za}}


%Telephone: (800) 555--1212, Fax: (888) 555--1212}


\maketitle

\begin{abstract}
Determining the optimal size and orientation of small-scale residential based PV arrays will become increasingly complex in the future smart grid environment with the introduction of smart meters and dynamic tariffs. However consumers can leverage the availability of smart meter data to conduct a more detailed exploration of PV investment options for their particular circumstances. In this paper, an optimization method for PV orientation and sizing is proposed whereby maximizing the PV investment value is set as the defining objective. Solar insolation and PV array models are described to form the basis of the PV array optimization strategy. A constrained particle swarm optimization algorithm is selected due to its strong performance in non-linear applications. The optimization algorithm is applied to real-world metered data to quantify the possible investment value of a PV installation under different energy retailers and tariff structures. The arrangement with the highest value is determined to enable prospective small-scale PV investors to select the most cost-effective system.
\end{abstract}


% \fi
\IEEEpeerreviewmaketitle

\section{Introduction}
Solar power generation, especially Photovoltaics (PV), has undergone enormous growth over the last decade. The integration of smart grid technologies, particularly smart metering, has the potential to significantly influence the penetration of small-scale rooftop PV systems in the distribution network.

Smart grid trials around the world suggest that the complexity of the business case for small-scale rooftop PV systems will increase in a smart grid environment. The Australian Government initiated `Smart Grid, Smart City' (SGSC) trial was one of the largest and widest ranging smart grid technology assessments to have been conducted in the world to date \cite{Arup14}. The SGSC study, and research conducted by other organizations such as the Grattan Institute \cite{Wood15}, found that under existing incentive schemes and flat electricity tariffs, customers are incentivized to install larger PV systems; with the remaining customers effectively providing a cross-subsidy for electricity costs of PV owners. In a smart grid environment, with smart meters and dynamic electricity tariffs, it was found that PV would continue to grow in Australia; however, the average system size would reduce for new systems. Consequently there will be a need for comprehensive and reasoned decision making tools to ensure PV is integrated to its maximum potential.

From a cumulative energy perspective, installation optimization has been extensively researched in literature \cite{Yadav13,Khatib16}. However existing optimization problems, such as the analysis undertaken by Koo et al. \cite{Koo16} and Rhodes et al. \cite{Rhodes14}, were found to be aimed largely towards determining the energetically optimal orientation rather than establishing an economic assessment based on temporal energy usage profiles. An economic assessment based on energy consumption was undertaken by Mulder et al. whereby the optimal size was investigated to leverage available feed-in tariffs and incentive schemes \cite{Mulder13}. However the high temporal resolution of smart meter data was not leveraged in the assessment; rather cumulative annual energy consumption data was utilized. 

In this paper, an exploration of PV installation optimization is conducted whereby the maximization of PV investment value is set as the defining objective and achieved through consideration of temporal energy usage (using smart meter data) in addition to other influential factors such as location, insolation data and electricity tariff structures. When consumers have a choice between flat tariff and dynamic time-of-use (TOU) rates from a multitude of energy retailers, the best plan is not self-evident, requiring an optimization strategy to select the most cost effective arrangement. This research aims to enable prospective small-scale PV investors to determine the highest value system and energy plan for their individual energy usage profile and circumstances.

In this paper, appropriate solar insolation and PV array models are presented to establish the underlying model for the optimization problem. Particle swarm optimization (PSO) is selected due to its speed and relative ease of application to non-linear problems. The basic PSO algorithm is modified with a penalty function to enable the handling of constraints related to feasible system size and position. Real-world solar insolation, electricity smart meter data and currently available Australian electricity plans from multiple retailers are used to quantify the investment value of an optimally sized and positioned PV array. Although the development and demonstration of the optimization problem is presented from an Australian energy market perspective, the methodology adopted in this paper is easily adaptable for other locations and countries.
\vspace*{+1.75\baselineskip} 

\section{Energy Generation Models}

Critical to PV system optimization is the definition of models for solar insolation and PV generated energy. The models considered in this paper are subsequently presented in this section.

\subsection{Solar Insolation Model}
\label{ss:solar}
While solar insolation with a high temporal resolution (e.g. 15 min, 30 min, 1 hour) is available for many locations, either through ground or satellite based measurements, such data is only available on a paid basis through software such as Meteonorm or SolarGIS. Freely available insolation databases, such as the Australian Bureau of Meteorology Climate Data Online (CDO) database \cite{BOM16}, typically only maintain daily global insolation. Consequently, when investigating PV array optimization based on hourly consumption data, a model for hourly insolation from daily data is required.

Irradiance incident on the horizontal plane can be divided into two components: beam (also known as direct) and diffuse. Consequently, the daily global insolation must first be separated into these components. Erbs et al. developed a correlation model for the fraction of daily diffuse insolation $H_d$, and daily global insolation $H$ for which the reader is referred to \cite{Erbs82} and \cite{Duffie13}. Determining an estimate of daily beam insolation, $H_b$, is then a simple process since

\begin{equation}\label{eq:dailyinsol}
H=H_{b}+H_{d}
\end{equation}

Collares-Pereira and Rabl \cite{Collares79}, established an estimate of the ratio of the global hourly insolation, $I$, and the daily global insolation as described by (\ref{eq:hourGHI}). In addition, Liu and Jordan described a similar estimate for the hourly diffuse insolation $I_d$ in \cite{Liu60} as shown in (\ref{eq:hourDHI}):

\begin{equation}\label{eq:hourGHI}
r_{t}=\dfrac{I}{H}=\dfrac{\pi}{24}(a+b\cos\omega)\dfrac{\cos\omega-\cos\omega_{s}}{\sin\omega-\dfrac{\pi\omega_{s}}{180}\cos\omega_{s}}
\end{equation}

\begin{equation}\label{eq:hourDHI}
r_{d}=\dfrac{I_{d}}{H_{d}}=\dfrac{\pi}{24}\dfrac{\cos\omega-\cos\omega_{s}}{\sin\omega-\dfrac{\pi\omega_{s}}{180}\cos\omega_{s}}
\end{equation}
where $a$ and $b$ are coefficients described by an empirical relationship, $\omega$ is the solar hour angle and $\omega_s$ is the sunset hour angle. Applying the daily estimated beam and diffuse insolation to (\ref{eq:hourGHI}) and (\ref{eq:hourDHI}), enables the hourly global and diffuse insolation to be determined. The beam hourly insolation can then be estimated by an equivalent hourly insolation version of (\ref{eq:dailyinsol}).

After establishing a model for the horizontal hourly insolation, a model of the insolation on the tilted plane must be developed. The total tilted plane hourly insolation ($I_T$), consists of beam ($I_b$), diffuse ($I_d$) and reflected insolation. Noorian et al. present an evaluation of twelve models for estimating hourly diffuse insolation on a tilted plane \cite{Noorian08}. The evaluation found that in general the Reindl et al. \cite{Reindl90} model, also known as the Hay-Davies-Klucher-Reindl (HDKR) model in \cite{Duffie13}, is one of the more accurate models and is relatively simple to use. Consequently, the HDKR model, described by (\ref{eq:tiltinsol}), was selected for this research.

\begin{IEEEeqnarray}{Rl}
I_{T}=&\left(I_{b}+A_{i}I_{d}\right)R_{b}\nonumber\\
&+\:(1-A_{i})\left(\frac{1+\cos\beta}{2}\right)\left[1+f\sin^{3}\left(\frac{\beta}{2}\right)\right]\nonumber\\
&+\:I\rho_{g}\left(\frac{1-\cos\beta}{2}\right)\label{eq:tiltinsol}%
\end{IEEEeqnarray}
where $A_i=I_b/I_o$, $f=\sqrt{I_b/I}$, $I_o$ is the hourly extraterrestrial insolation that would be incident on a horizontal plane projected from the Earth's surface and $\rho_g$ is the ground reflectance, assumed to be 0.2 in this research representing short grassland. The ratio of tilted to horizontal beam radiation ratio $R_b$ is defined as

\begin{equation}\label{eq:beamratio}
R_{b}=\frac{\cos\theta}{\cos\theta_{z}}
\end{equation}
where
\begin{equation}\label{eq:aoi}
\begin{split}
\begin{aligned}
\cos\theta= & \sin\delta\sin\phi\cos\beta+\sin\delta\cos\phi\sin\beta\cos\gamma\\
 & +\cos\delta\cos\phi\cos\beta\cos\omega\\
 & -\cos\delta\sin\phi\sin\beta\cos\gamma\cos\omega\\
 & +\cos\delta\sin\beta\sin\gamma\sin\omega
\end{aligned}
\end{split}
\end{equation}

\begin{equation}\label{eq:zenith}
\cos\theta_{z}=\cos\phi\cos\delta\cos\omega+\sin\phi\sin\delta
\end{equation}

The angles $\theta$, $\theta_z$, $\delta$, $\phi$, $\beta$ and $\gamma$ refer to the beam irradiance angle of incidence, zenith angle of the sun, solar declination, latitude, panel tilt and panel azimuth respectively. Note that (\ref{eq:aoi}) has been developed for locations in the southern hemisphere whereby $\gamma = 0$ implies north facing surfaces.

Equation (\ref{eq:tiltinsol}) constitutes the defining insolation equation for the optimization problem.

\subsection{Photovoltaic Array Model}
\label{ss:photo}

The energy balance equation for a PV module, cooled by losses to the surroundings, is described by (\ref{eq:pvcool}) \cite{Duffie13}:

\begin{equation}\label{eq:pvcool}
(\tau\alpha)G_{T}=\eta_{c}G_{T}+U_{L}(T_{c}-T_{a})
\end{equation}
where ($\tau\alpha$) is the transmittance-absorbance product, $\eta_c$ is the module efficiency quantifying the effectiveness of converting irradiance to electrical energy, $T_c$  is the cell temperature, $T_a$ is the ambient temperature, $G_T$ is the incident irradiance and $U_L$ is the heat loss coefficient.

Manufacturer datasheets detail a PV module's performance under both standard test conditions (STC) and nominal operating cell temperature (NOCT) conditions. The NOCT conditions include incident irradiance of 800 W/m$^2$ and an ambient temperature of 20$^\circ$C with the PV cells under no load. Following the steps described in \cite{Duffie13}, from (\ref{eq:pvcool}) and through knowledge of the STC and NOCT performance, the temperature $T_c$ of a PV panel can be defined as: 

\begin{equation}\label{eq:celltemp}
T_{c}=T_{a}+(T_{NOCT}-20)\cdot\frac{G_{T}}{800}\cdot(1-\eta_{mpp,STC})
\end{equation}
where $T_{NOCT}$ and $\eta_{mpp,STC}$ are the temperature of the panel and panel efficiency under NOCT and STC conditions respectively.

Assuming $G_T$ is constant, the output energy of a PV array over an hour period is described by:

\begin{equation}\label{eq:pvpower}
E_{pv}=A_{c}ZI_{T}\eta_{mpp}\eta_e
\end{equation}
where $Z$ in the number of panels, $A_c$ is the panel area and $\eta_e$ is the efficiency of the associated balance of plant. The efficiency $\eta_{mpp}$ of the PV array at particular operating and environmental conditions is defined as:

\begin{equation}\label{eq:opeff}
\eta_{mpp}=\eta_{mpp,STC}+\mu_{mpp}(T_{c}-T_{a})
\end{equation}
where $\mu_{mpp}$ (\%/W) is the power coefficient detailed on the manufacturer's datasheet. 

Equations~(\ref{eq:pvpower}) and (\ref{eq:opeff}) constitute the defining PV array model to be used in the objective function of the optimization problem described in the next section.

\section{Optimization Problem}
\label{s:optim}

The objective of the optimization problem, defined in the following subsection, is to maximize the value of the PV investment option by selecting the optimal size and position of a proposed PV array given a known lowest cost `do nothing' electricity tariff. Through a comparison of optimized systems for different electricity tariff structures and energy retailers, the optimal investment option can be determined. Due to the relatively long investment period, the time value of money must be taken into consideration through a net present value (NPV) analysis.

\subsection{Problem Definition}
\label{ss:probdef}

\vspace*{+0.5\baselineskip}
\noindent\emph{Given:}
\begin{enumerate}
\item Maximum allowable PV panels ($Z_{max}=30$)
\item PV cost per watt peak ($U_{pv} = \$2.30$/W$_{p}$)
\item Latitude and longitude of the location
\item Hourly load and insolation profile
\item A real discount rate of 6\%
\item An inflation rate of 2\% 
\item PV system balance of plant efficiency (90\%)
\item System lifespan (15 years)
\end{enumerate}

\vspace*{+0.5\baselineskip}
\noindent\emph{Find:} Tilt angle $\beta$, azimuth angle $\gamma$ and number of panels $Z$

\vspace*{+0.5\baselineskip}
\noindent\emph{Objective:}

\begin{IEEEeqnarray}{Rl}
\max_{\beta,\gamma,Z} NPV=&\sum_{q=1}^{Q}\frac{\bigl(C_{base,q}-C_{pv,q}(\beta,\gamma,Z)\bigr)\left(1+r_{i,eff}\right)^q}{\left(1+r_{d,eff}\right)^{q}}\nonumber\\
&-\:S_{pv}(Z)\label{eq:objfunc}%
\end{IEEEeqnarray}

\noindent\emph{Subject to:}
\begin{IEEEeqnarray}{rCll}
g_{1}(\beta)&=&\vert\beta-90\vert-90 \leq 0&\text{for}\;\beta\in\mathbb{R}\IEEEyesnumber\IEEEyessubnumber\label{eq:tiltcons}\\
g_{2}(\gamma)&=&\vert\gamma\vert-180 \leq 0&\text{for}\;\gamma\in\mathbb{R}\IEEEyessubnumber\label{eq:azicons}\\
g_{3}(Z)&=&\vert Z-Z_{max}\vert-Z_{max} \leq 0\;&\text{for}\;Z\in\mathbb{Z}^{+}\IEEEyessubnumber\label{eq:panelcons}%
\end{IEEEeqnarray}

In (\ref{eq:objfunc}), $C_{base,q}$ is the cost of energy without PV and $C_{pv,q}$ is the cost of electricity with PV within the billing period $q$ (the difference of which constitutes the monetary savings achieved through the installation of PV). $S_{pv}$ and $Q$ are the PV system cost and the number of billing periods (assume quarterly billing with 60 quarters over 15 years). Given the quarterly billing cycle, the annual inflation and discount rates of 2\% and 6\% respectively must be are adjusted to the quarterly effective rates $r_{i,eff}$ and $r_{d,eff}$. 

The terms $C_{base,q}$, $C_{pv,q}$ and $S_{pv}$ are defined by equations (\ref{eq:basecost}), (\ref{eq:pvcost}) and (\ref{eq:systemcost}) respectively.

\begin{equation}\label{eq:basecost}
C_{base,q}=\sum_{d=1}^{D}\left[ \sum_{h=1}^{24}\left( T_{grid0,h}E_{load,h}\right) +T_{sc0,d}\right]
\end{equation}

\begin{IEEEeqnarray}{Rl}
C_{pv,q}=&\sum_{d=1}^{D}\Bigg\lbrace T_{sc,d}+\sum_{h=1}^{H}\Big[T_{grid,h}\max\bigl(0,E_{bal,h}(\beta,\gamma,Z)\bigr)\nonumber\\
&-\;T_{feed,h}\max\bigl(0,-E_{bal,h}(\beta,\gamma,Z)\bigr)\Big]\Bigg\rbrace\label{eq:pvcost}
\end{IEEEeqnarray}

\begin{equation}\label{eq:energybal}
E_{bal,h}(\beta,\gamma,Z)=E_{load,h}-E_{pv,h}(\beta,\gamma,Z)
\end{equation}
\begin{equation}\label{eq:systemcost}
S_{pv}=U_{pv}P_{pv,rat}(Z)-M_{loc}P_{pv,rat}(Z)C_{STC}
\end{equation}

In~(\ref{eq:basecost}) and (\ref{eq:pvcost}), $T_{grid0,h}$ and $T_{grid,h}$ are the grid imported electricity tariff of the base plan and tested plan respectively for the $h^{th}$ hour of day $d$, with $D$ days in the billing period. $T_{sc0,d}$ and $T_{sc,d}$ are the daily electricity supply charges for the base plan and tested plan respectively. $E_{bal,h}$ is the net energy flow balance defined in (\ref{eq:energybal}). $E_{load,h}$, $E_{pv,h}$ and $T_{feed,h}$ are the energy consumed by the load, the energy generated by the PV system and the feed-in tariff respectively. 

In an Australian context, the total PV system cost $S_{pv}$ of rated power $P_{pv,rat}$ is reduced through an effective rebate provided through the small-scale technology scheme. Under this scheme, small-scale technology certificates (STCs) are generated based on the system size and a location multiplier $M_{loc}$ (assumed to be 20.73 in this paper). It is assumed the STCs are worth $C_{STC}$ = \$35 to the system owner.

\begin{figure}[!t]
\centering
\includegraphics[width=\linewidth]{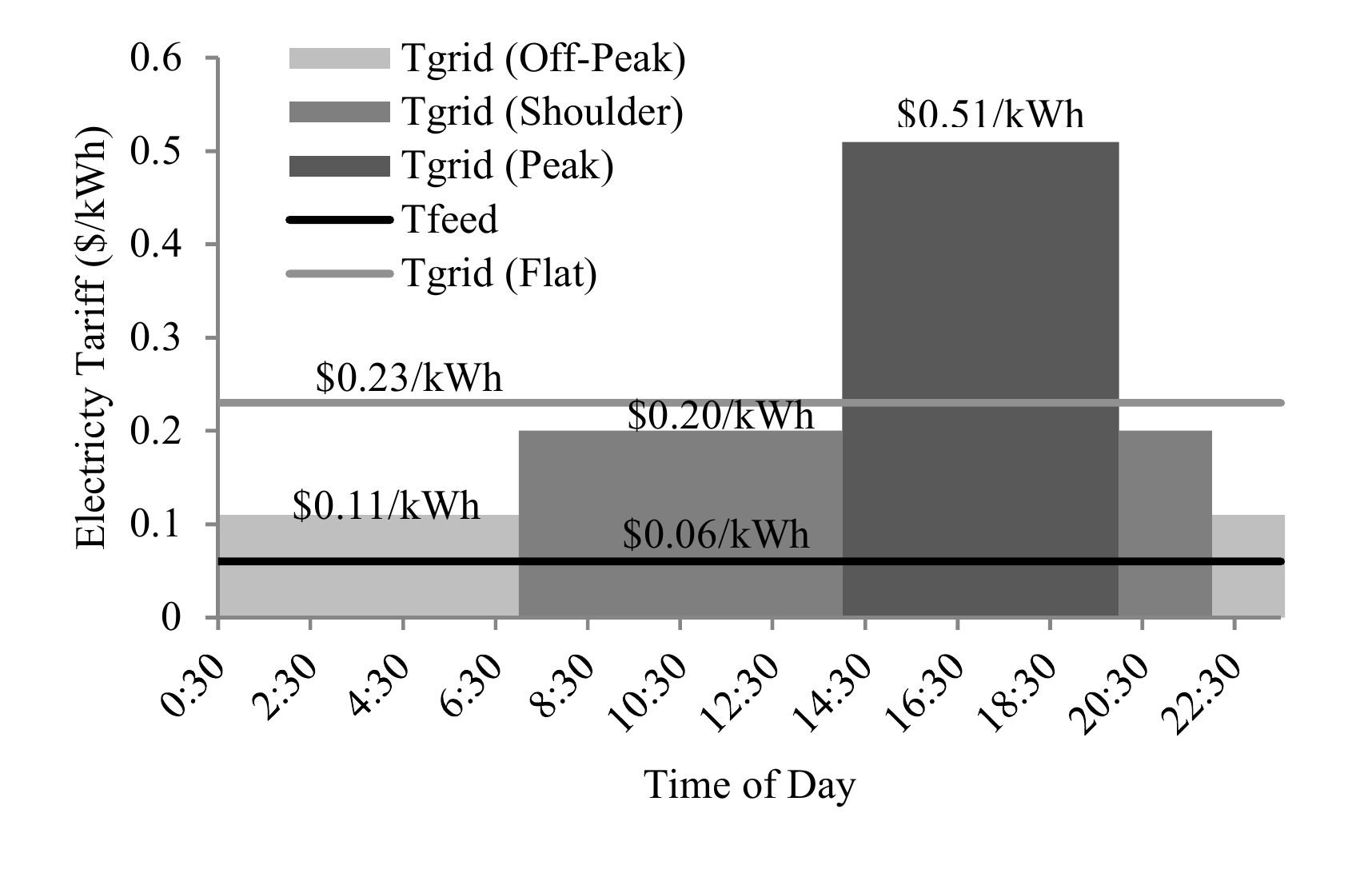}
\caption{Typical electricity tariff structure and rates}
\label{fig:tariff}
\end{figure}

\subsection{Optimization Algorithm}
\label{ss:optalg}

Due to its speed, global search performance and relative simplicity of application \cite{Sun12}, PSO has seen many applications in PV optimization problems \cite{Yadav13}. Particle swarm optimization (PSO) is a metaheuristic programming method simulating the social interaction within bird flocks and fish schools to achieve a global objective in the absence of centralized control. Each individual swarm agent is identified as a particle that flies through the solution space, defined by position and velocity vectors. The dimensionality, $J$, of the vectors is equal to the number of optimization parameters (for the problem defined in this paper, $J=3$). Consequently for the $i^{th}$ particle in the $n^{th}$ iteration, the position $\boldsymbol{x}_{i,n}$ and velocity $\boldsymbol{v}_{i,n}$ vectors are defined as:

\begin{equation}\label{eq:position}
\boldsymbol{x}_{i,n}=(x_{i,n}^1,x_{i,n}^2,x_{i,n}^3)=(\beta_{i,n},\gamma_{i,n},Z_{i,n})
\end{equation}

\begin{equation}\label{eq:velocity}
\boldsymbol{v}_{i,n}=(v_{i,n}^1,v_{i,n}^2,v_{i,n}^3)
\end{equation}

In PSO, the position $x_{i,n}^j$ and velocity $v_{i,n}^j$ of each particle are updated component-wise in iteration $n+1$ through knowledge of the particle's personal best position and the global best position within the swarm as expressed in (\ref{eq:velupdate}) and (\ref{eq:posupdate}):

\begin{IEEEeqnarray}{Rl}
v_{i,n+1}^j=&\chi\Big[v_{i,n}^j+c_1r_{i,n}^j\bigl(P_{i,n}^j-x_{i,n}^j\bigr)\nonumber\\
&+\:c_2R_{i,n}^j\bigl(G_n^j-x_{i,n}^j\bigr)\Big]\label{eq:velupdate}
%\mzreset
\end{IEEEeqnarray}

\begin{equation}\label{eq:posupdate}
x_{i,n+1}^j=x_{i,n}^j+v_{i,n+1}^j   
\end{equation}
where, for dimension $j$ of particle $i$, $P_{i,n}^j$ is the personal best position, $G_{n}^j$ is the global best position of the swarm, $c_1$ and $c_2$ are acceleration coefficients and $r_{i,n}^j$ and $R_{i,n}^j$ are two sequences of random numbers distributed uniformly over $(0,1)$. Equation~(\ref{eq:velupdate}) represents a modification of the basic PSO algorithm by inclusion of the term $\chi$. Known as PSO with constriction factor (PSO-Co), the modification gives significantly improved convergence performance over the basic form \cite{Sun12}. In this research, $c_1$, $c_2$ and $\chi$ are set be 2.05, 2.05 and 0.729 respectively according to the recommendations in \cite{Sun12}. The algorithm is terminated when the maximum number of iterations $N$ is reached or the global best position is sufficiently close to the actual solution.

In order to handle optimization constraints, the most common method is the introduction of a penalty function to the underlying objective function \cite{Sun12}. A penalty function of the form explored by Parsopoulos and Vrahatis \cite{Parsopoulos02} is considered in this research in order to handle the PV array size and position constraints. The objective function is now considered to take the following form:

\begin{equation}\label{eq:objfuncpen}
F(\boldsymbol{x})=f(\boldsymbol{x})+H(\boldsymbol{x},n)
\end{equation}
where $\boldsymbol{x}$ is the optimization variables vector such that $\boldsymbol{x}=(\beta,\gamma,Z)$, $f(\boldsymbol{x})$ is the original objective function defined by (\ref{eq:objfunc}) and $H(\boldsymbol{x})$ is a penalty function of the form:

\begin{equation}\label{eq:penalty}
H(\boldsymbol{x},n)=h(n)\sum_{j=1}^{J}\theta\bigl(y_{j}(\boldsymbol{x})\bigr)y_j(\boldsymbol{x})^{\alpha\bigl(y_j(\boldsymbol{x})\bigr)}
\end{equation}

In (\ref{eq:penalty}), $y_j(\boldsymbol{x})=\max\bigl(0,g_j(\boldsymbol{x})\bigr)$ is a relative violated function of the constraints $g_j(\boldsymbol{x})$ in (\ref{eq:tiltcons})\textendash(\ref{eq:panelcons}) (such that the penalty function is zero when all optimization variables satisfy the constraints); $\theta\bigl( y_j(\boldsymbol{x})\bigr)$ is a multi-stage assignment function (which scales the penalty depending on the value of $y_j(\boldsymbol{x})$); $\alpha\bigl(y_j(\boldsymbol{x})\bigr)$ is the power of the penalty function; and $h(n)$ is a dynamically modified weight factor. The penalty parameters in (\ref{eq:penalty}) are problem dependent however the values defined by Parsopoulos and Vrahatis in \cite{Parsopoulos02} were found to be suitable.

As the optimization parameter $Z$ is limited to integer values while $\beta$ and $\gamma$ may take any real value within the domain of the constraints, the problem is classified as a mixed-integer non-linear programming problem. In order to handle the discrete parameter $Z$, the hypercube nearest-vertex approach adopted in \cite{Chowdhury13} is utilized. For installation simplicity, the tilt and azimuth angles were also considered as discrete in this analysis. 

The PSO optimization algorithm was developed and simulated in Matlab version R2015b.

\section{Input Data}
\label{s:input}

Between 2012\textendash2014 the SGSC project collected smart meter energy data for approximately 13,700 residences in the Sydney, Newcastle and Hunter Valley regions of New South Wales, Australia \cite{Arup14}. Electricity consumption data measured over a one-year period from three arbitrarily selected customers within the SGSC database were used to demonstrate the optimization strategy developed in this paper.

Five years of daily insolation and ambient temperature from the nearest weather stations within the CDO database were used in the analysis. Daily maximum temperature data was used due to the lack of hourly data, yielding conservative estimates for PV performance.

The electricity tariff structures tested in the optimization problem were based on real 2016 rates for three large Australian retailers. For each retailer, a flat tariff and a TOU tariff were considered for which example rates are shown in Fig.~\ref{fig:tariff}.

The PV arrays were modelled based on 250 W Trina TSM-PD05.05 polycrystalline PV modules.

\section{Results and Discussion}
\label{s:results}

\begin{table*}[!b]
\renewcommand{\arraystretch}{1.3}
\caption{Characteristics and Economic Performance of Optimized PV Systems under Different Retail Electricity Plans}
\label{tab:results}
\centering
\begin{tabular}{|c|c|c|c|c|c|c|c|c|c|}
\hline
Customer & Retailer & Tariff & Size (kW) & Tilt & Azimuth & NPV  & MIRR & Payback (Years) & Plan Saving \\\hline\hline
\multirow{6}{*}{1} & A & TOU & 5.76 & $28^\circ$ & $32^\circ$ & \$4,861 & 6.29\% & 9.3 & \$1,073 \\\hhline{~*{9}{-}}
 & A & Flat & 4.51 & $26^\circ$ & $8^\circ$ & \$4,965 & 7.00\% & 8.3 & \$1,176 \\\hhline{~*{9}{-}}
  & \cellcolor[gray]{0.8}B & \cellcolor[gray]{0.8}TOU & \cellcolor[gray]{0.8}4.26 & \cellcolor[gray]{0.8}$29^\circ$ & \cellcolor[gray]{0.8}$35^\circ$ & \cellcolor[gray]{0.8}\$5,606 & \cellcolor[gray]{0.8}7.56\% & \cellcolor[gray]{0.8}7.8 & \cellcolor[gray]{0.8}\$1,817 \\\hhline{~*{9}{-}}
 & B & Flat & 4.01 & $25^\circ$ & $9^\circ$ & \$3,983 & 6.71\% & 8.8 & \$194 \\\hhline{~*{9}{-}}
 & C & TOU  & 4.51 & $29^\circ$ & $37^\circ$ & \$4,850 & 6.94\% & 8.3 & \$1,062 \\\cline{2-10}
 & C & Flat & 4.01 & $24^\circ$ & $10^\circ$ & \$3,789 & 6.58\% & 9.0 & \$0 \\\hline\hline
 
\multirow{6}{*}{2} & A & TOU & - & - & - & - & - & - & - \\\hhline{~*{9}{-}}
 & \cellcolor[gray]{0.8}A & \cellcolor[gray]{0.8}Flat  & \cellcolor[gray]{0.8}1.75 & \cellcolor[gray]{0.8}$28^\circ$ & \cellcolor[gray]{0.8}$4^\circ$ & \cellcolor[gray]{0.8}\$1,431 & \cellcolor[gray]{0.8}6.20\% & \cellcolor[gray]{0.8}9.4 & \cellcolor[gray]{0.8}\$1,431 \\\hhline{~*{9}{-}}
 & B & TOU & 1.75 & $30^\circ$ & $29^\circ$ & \$404 & 4.24\% & 13.2 & \$404 \\\hhline{~*{9}{-}}
 & B & Flat & 1.50 & $28^\circ$ & $3^\circ$ & \$1,389 & 6.51\% & 9.1 & \$1,389 \\\hhline{~*{9}{-}}
 & C & TOU  & - & - & - & - & - & - & - \\\hhline{~*{9}{-}}
 & C & Flat & 1.50 & $28^\circ$ & $3^\circ$ & \$1,204 & 6.15\% & 9.7 & \$1,204 \\\hline\hline  

\multirow{6}{*}{3} & A & TOU  & 2.76 & $27^\circ$ & $33^\circ$ & \$2,430 & 6.40\% & 9.3 & \$30 \\\hhline{~*{9}{-}}
 & \cellcolor[gray]{0.8}A & \cellcolor[gray]{0.8}Flat & \cellcolor[gray]{0.8}2.26 & \cellcolor[gray]{0.8}$25^\circ$ & \cellcolor[gray]{0.8}$3^\circ$  & \cellcolor[gray]{0.8}\$2,853 & \cellcolor[gray]{0.8}7.43\% & \cellcolor[gray]{0.8}7.8 & \cellcolor[gray]{0.8}\$453 \\\hhline{~*{9}{-}}
 & B & TOU  & 2.56 & $28^\circ$ & $35^\circ$ & \$2,778 & 7.34\% & 8.1 & \$379 \\\hhline{~*{9}{-}}
 & B & Flat & 2.00 & $25^\circ$ & $3^\circ$  & \$2,712 & 7.65\% & 7.6 & \$312 \\\hhline{~*{9}{-}}
 & C & TOU  & 2.26 & $28^\circ$ & $37^\circ$ & \$2,399 & 6.91\% & 8.6 & \$0 \\\hhline{~*{9}{-}}
 & C & Flat & 2.00 & $25^\circ$ & $3^\circ$  & \$2,544 & 7.44\% & 7.8 & \$145 \\\hline
\end{tabular}
\end{table*}

Table~\ref{tab:results} summarizes the results of the optimization problem of Section~\ref{ss:probdef} under both flat and TOU tariffs from the three different Australian energy retailers considered. For each customer, installing PV was found to be a positive investment with at least one of the available electricity retail plans providing a positive NPV. The retail plans yielding the highest NPV for each customer are highlighted in gray. 

Referring to Table~\ref{tab:results}, Customer~2 and Customer~3 were found to achieve the greatest benefit from a flat tariff retail plan from Retailer~A. However Customer~1 was found to benefit most from a TOU tariff from Retailer~B. The diversity in NPV achievable amongst the customers and the corresponding optimal PV system sizes is also evident. For Customer~1, as the beneficiary of the largest potential benefit, the maximum NPV was achieved with a system size of 4.26 kW. In comparison, Customer~2 and Customer~3 benefit from relatively smaller sizes of 1.75 kW and 2.26 kW respectively.

Evident in the final column of Table~\ref{tab:results}, significant savings can be achieved by considering available retail electricity plan options in the analysis. In the hypothetical situation whereby Customer~1 previously held a flat plan from Retailer~C prior to the decision to invest in PV, optimizing the PV system for this particular plan would yield a NPV of \$3,789 from a 4.01 kW system. However when alternative retail plans are considered in the analysis, a far higher NPV of \$5,606 is achievable, representing an additional saving of \$1,817 against the previous plan. It should be noted that under TOU tariffs from Retailer~A and Retailer~C, no PV system could be found that resulted in a net benefit for Customer~2. Consequently the entries in Table~\ref{tab:results} were left empty.

Although positive NPV was set as the objective of the optimization problem, a positive NPV should not be considered alone in the investment analysis. Consequently common economic metrics of modified internal rate of return (MIRR) and payback period were also calculated and summarized in Table~\ref{tab:results}. In order for an investment in PV to be considered preferable to other investment alternatives, the MIRR should be higher than the discount rate. As can be clearly seen in Table~\ref{tab:results}, the MIRRs associated with the systems yielding the largest NPV for each customer are greater than the assumed 6\% discount rate, with up to 7.56\% observed for the optimal system of Customer~1. Consequently for each customer, an optimally selected PV system and associated least-cost electricity plan was found to be a greater investment opportunity than the other market alternatives.

Payback periods well within the assumed system lifespan of 15 years were found for each customer. The optimal systems for Customer 1 and Customer 3 were found to have a payback period of 7.8 years while Customer 2 would experience a slightly longer period of 9.4 years.

\begin{figure}[!t]
\centering
\includegraphics[width=\linewidth]{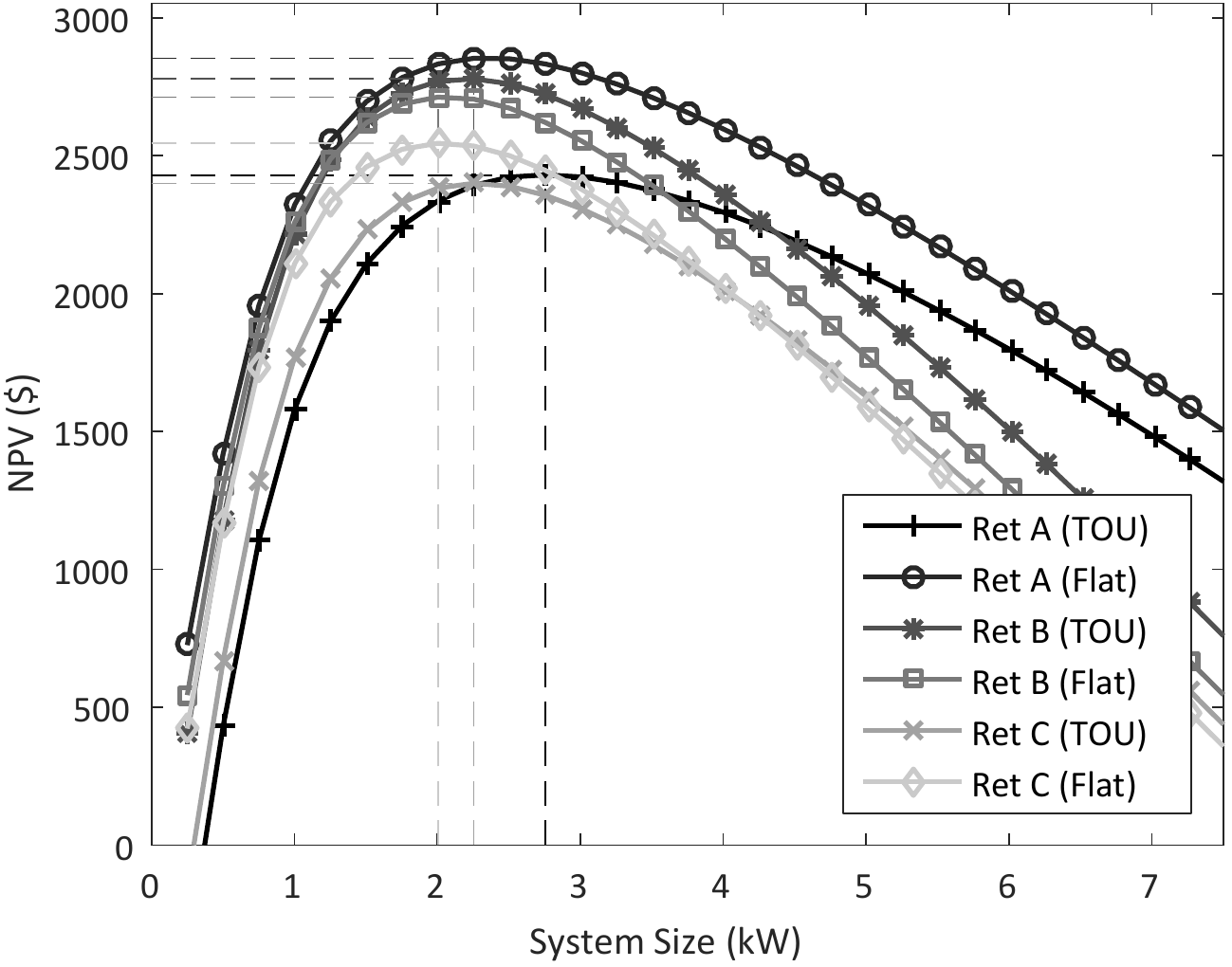}
\caption{NPV sensitivity to system size (Customer 3)}
\label{fig:ori_sensi}
\end{figure}

\begin{figure}[!t]
\centering
\includegraphics[width=\linewidth]{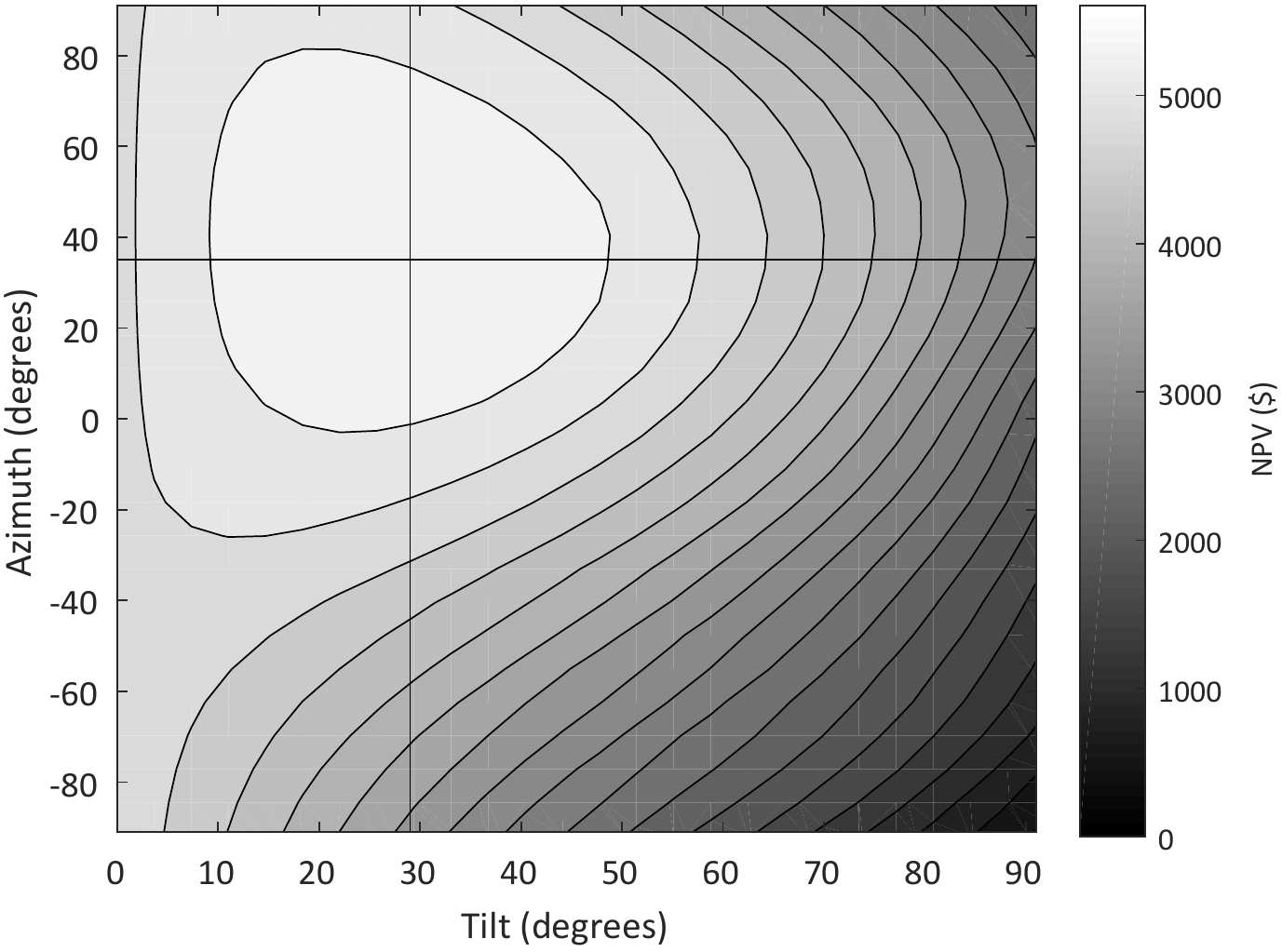}
\caption{NPV sensitivity to tilt and azimuth (Customer 1)}
\label{fig:size_sensi}
\end{figure}

The effect of sub-optimal system sizing and orientation is shown in Fig.~\ref{fig:size_sensi} and Fig.~\ref{fig:ori_sensi} respectively. For Customer~3, the sensitivity of NPV to different system sizes is clearly evident in Fig.~\ref{fig:size_sensi}. As the size is increased from one PV panel to the optimal size of 2.26 kW, the potential NPV rapidly increases. However if a larger system is selected, the NPV begins to decrease. Given the increased investment cost of a larger system, the rate of return and payback periods for the larger system would also be detrimentally effected.

The contour plot of Fig.~\ref{fig:ori_sensi}, shows an optimal position significantly skewed west of north (positive azimuth) for Customer~1. However it can be seen that NPV is relatively insensitive to minor deviations form the the sub-optimal orientation. If a north facing azimuth (0$^\circ$) and a tilt anywhere between 10$^\circ$\textendash40$^\circ$ are selected, the resulting NPV would be relatively close to the NPV of the optimal system. Therefore it may be concluded that the two primary considerations for economic optimization are the system size and accompanying energy plan. 

The diversity in the optimal system sizes and orientations amongst both the customers and available retail plans, highlights the necessity to adopt an optimization strategy, such as the one presented in this paper, during the investment decision process. 

\section{Conclusion}

The lack of transparency regarding the highest value PV system size, orientation and retail electricity plans for a given location, insolation and power consumption profile is a key driver for the development of a PV optimization strategy. 

In this paper, hourly insolation and PV array models are presented and incorporated in an objective function aiming to maximize the net present value of a PV investment by optimally sizing and positioning the array and through selection of the least-cost electricity plan amongst a range of alternative options. A modified particle swarm optimization algorithm is utilized to solve the optimization problem.

For each of the customers assessed, an optimal PV system producing a positive NPV (up to \$5,606) and yielding a rate of return greater than the cost of capital was found (up to 7.56\%). Furthermore, the importance of least-cost energy plan selection was established with up to \$1,817 in additional savings found.

The diversity of PV system sizes and potential net benefits observed amongst the customers and energy plans assessed, highlights the necessity for an optimization tool prior to selection of a residential PV system. To establish a wider assessment of economic performance metrics under the current Australian market conditions, a wider sample of customer data is necessary. Such an assessment is currently the subject of further research by the authors.


%\section*{Acknowledgment}



%\bibliographystyle{IEEEtran}

\bibliographystyle{IEEEtran}
\bibliography{IEEEabrv,icrera2016}

\begin{thebibliography}{10}
\providecommand{\url}[1]{#1}
\csname url@samestyle\endcsname
\providecommand{\newblock}{\relax}
\providecommand{\bibinfo}[2]{#2}
\providecommand{\BIBentrySTDinterwordspacing}{\spaceskip=0pt\relax}
\providecommand{\BIBentryALTinterwordstretchfactor}{4}
\providecommand{\BIBentryALTinterwordspacing}{\spaceskip=\fontdimen2\font plus
\BIBentryALTinterwordstretchfactor\fontdimen3\font minus
  \fontdimen4\font\relax}
\providecommand{\BIBforeignlanguage}[2]{{%
\expandafter\ifx\csname l@#1\endcsname\relax
\typeout{** WARNING: IEEEtran.bst: No hyphenation pattern has been}%
\typeout{** loaded for the language `#1'. Using the pattern for}%
\typeout{** the default language instead.}%
\else
\language=\csname l@#1\endcsname
\fi
#2}}
\providecommand{\BIBdecl}{\relax}
\BIBdecl

\bibitem{Arup14}
M.~Norris, R.~Cliff, R.~Sharp, S.~Koci, and H.~Gardner, ``{Smart Grid, Smart
  City}: Shaping {Australia's} energy future national cost benefit
  assessment,'' Arup, Report, July 2014.

\bibitem{Wood15}
T.~Wood, D.~Blowers, and C.~Chisholm, ``Sundown, sunrise: How {Australia} can
  finally get solar power right,'' Grattan Institute, Report 2015-2, May 2015.

\bibitem{Yadav13}
A.~K. Yadav and S.~S. Chandel, ``Tilt angle optimization to maximize incident
  solar radiation: A review,'' \emph{Renewable and Sustainable Energy Reviews},
  vol.~23, pp. 503--513, 2013.

\bibitem{Khatib16}
T.~Khatib, I.~A. Ibrahim, and A.~Mohamed, ``A review on sizing methodologies of
  photovoltaic array and storage battery in a standalone photovoltaic system,''
  \emph{Energy Conversion and Management}, vol. 120, pp. 430--448, 2016.

\bibitem{Koo16}
C.~Koo, T.~Hong, M.~Lee, and J.~Kim, ``An integrated multi-objective
  optimization model for determining the optimal solution in implementing the
  rooftop photovoltaic system,'' \emph{Renewable and Sustainable Energy
  Reviews}, vol.~57, pp. 822--837, 2016.

\bibitem{Rhodes14}
J.~D. Rhodes, C.~R. Upshaw, W.~J. Cole, C.~L. Holcomb, and M.~E. Webber, ``A
  multi-objective assessment of the effect of solar{ PV} array orientation and
  tilt on energy production and system economics,'' \emph{Solar Energy}, vol.
  108, no.~0, pp. 28--40, 2014.

\bibitem{Mulder13}
G.~Mulder, D.~Six, B.~Claessens, T.~Broes, N.~Omar, and J.~V. Mierlo, ``The
  dimensioning of {PV}-battery systems depending on the incentive and selling
  price conditions,'' \emph{Applied Energy}, vol. 111, pp. 1126--1135, 2013.

\bibitem{BOM16}
\BIBentryALTinterwordspacing
{Australian Government Bureau of Meteorology}, ``Climate data online,'' 2016.
  [Online]. Available: \url{http://www.bom.gov.au/climate/data/}
\BIBentrySTDinterwordspacing

\bibitem{Erbs82}
D.~G. Erbs, S.~A. Klein, and J.~A. Duffie, ``Estimation of the diffuse
  radiation fraction for hourly, daily and monthly-average global radiation,''
  \emph{Solar Energy}, vol.~28, no.~4, pp. 293--302, 1982.

\bibitem{Duffie13}
J.~A. Duffie and W.~A. Beckman, \emph{Solar engineering of thermal
  processes}.\hskip 1em plus 0.5em minus 0.4em\relax Hoboken, N.J.: Wiley,
  2013, vol. 4th.

\bibitem{Collares79}
M.~Collares-Pereira and A.~Rabl, ``The average distribution of solar
  radiation-correlations between diffuse and hemispherical and between daily
  and hourly insolation values,'' \emph{Solar Energy}, vol.~22, no.~2, pp.
  155--164, 1979.

\bibitem{Liu60}
B.~Y.~H. Liu and R.~C. Jordan, ``The interrelationship and characteristic
  distribution of direct, diffuse and total solar radiation,'' \emph{Solar
  Energy}, vol.~4, no.~3, pp. 1--19, 1960.

\bibitem{Noorian08}
A.~M. Noorian, I.~Moradi, and G.~A. Kamali, ``Evaluation of 12 models to
  estimate hourly diffuse irradiation on inclined surfaces,'' \emph{Renewable
  Energy}, vol.~33, no.~6, pp. 1406--1412, 2008.

\bibitem{Reindl90}
D.~T. Reindl, W.~A. Beckman, and J.~A. Duffie, ``Evaluation of hourly tilted
  surface radiation models,'' \emph{Solar Energy}, vol.~45, no.~1, pp. 9--17,
  1990.

\bibitem{Sun12}
J.~Sun, C.-H. Lai, and X.-J. Wu, \emph{Particle swarm optimisation}, ser.
  Chapman and Hall CRC numerical analysis and scientific computing.\hskip 1em
  plus 0.5em minus 0.4em\relax Boca Raton, Fla.: CRC Press, 2012.

\bibitem{Parsopoulos02}
K.~E. Parsopoulos and M.~N. Vrahatis, ``Particle swarm optimization method for
  constrained optimization problems,'' \emph{Intelligent Technologies~- Theory
  and Application: New Trends in Intelligent Technologies}, vol.~76, pp.
  214--220, 2002.

\bibitem{Chowdhury13}
S.~Chowdhury, W.~Tong, A.~Messac, and J.~Zhang, ``A mixed-discrete particle
  swarm optimization algorithm with explicit diversity-preservation,''
  \emph{Structural and Multidisciplinary Optimization}, vol.~47, no.~3, pp.
  367--388, 2013.

\end{thebibliography}

\end{document}